\documentclass[final,onefignum]{siamltex}
\usepackage{amsfonts}
\usepackage{amssymb}
\usepackage{graphicx}
\usepackage{latexsym}
\usepackage{times}
\usepackage{subeqn}
\usepackage{setspace}
\usepackage[sort,square]{natbib}
\usepackage[pdfborderstyle={/S/U/W 0}]{hyperref}
\usepackage{hyperref}

\usepackage{caption}

\newcommand{\Frac}[2]{\displaystyle \frac{\displaystyle #1}{\displaystyle #2}}
\newcommand{\sech}{ {\rm sech} }

\def\ds{\displaystyle} 

\def\sss{\scriptscriptstyle}
 
\def\R{{\bf R}}

\newcommand{\ttt}{^{\sss\top}} 

\author{Alessandro Colombo\footnotemark[2]\ \footnotemark[3] \and Fabio Dercole\footnotemark[2]}

\title{Discontinuity induced bifurcations of non-hyperbolic cycles in nonsmooth systems}

\begin{document}

\maketitle

\renewcommand{\thefootnote}{\fnsymbol{footnote}}
\footnotetext[2]{DEI, Politecnico di Milano, Via Ponzio 34/5, 20133 Milano, Italy,   [alessandro.colombo]fabio.dercole@polimi.it}
\footnotetext[3]{To whom correspondence should be addressed, Ph: +39 02 2399 4034; Fax: +39 02 2399 3412}
\renewcommand{\thefootnote}{\arabic{footnote}}

\begin{abstract}
We analyse three codimension-two bifurcations occurring in nonsmooth systems, when a non-hyperbolic cycle (fold, flip, and Neimark-Sacker cases, both in continuous- and discrete-time) interacts with one of the discontinuity boundaries characterising the system's dynamics. Rather than aiming at a complete unfolding of the three cases, which would require specific assumptions on both the class of nonsmooth system and the geometry of the involved boundary, we concentrate on the geometric features that are common to all scenarios.  We show that, at a generic intersection between the smooth and discontinuity induced bifurcation curves, a third curve generically emanates tangentially to the former.  This is the discontinuity induced bifurcation curve of the secondary invariant set (the other cycle, the double-period cycle, or the torus, respectively) involved in the smooth bifurcation.  The result can be explained intuitively, but its validity is proved here rigorously under very general conditions.  Three examples from different fields of science and engineering are also reported.
\end{abstract}



\pagestyle{myheadings}
\thispagestyle{plain}
\markboth{A. Colombo and F. Dercole}{Discontinuity induced bifurcations of non-hyperbolic cycles}

\section{Introduction}
\label{sec:int}
This article deals with the analysis of three particular codimension-two bifurcations in nonsmooth systems. Broadly speaking, nonsmooth systems are continuous- or discrete-time dynamical systems featuring some kind of discontinuity in the right-hand side of their governing equations whenever the system's state reaches a {\it discontinuity boundary}.
More specifically, nonsmooth systems include several classes, e.g., piecewise smooth \citep{Filippov88, Bernardo08}, impacting \citep{Brogliato99}, and hybrid \citep{Branicky98b, Lygeros03} systems, which have been largely used in the last decades as models in various fields of science and engineering (see references above and therein).

While methods of numerical continuation allow to easily detect and trace
bifurcation curves in two-parameter planes, understanding the geometry
of bifurcation curves around codimension-two points is a key to the
construction of complex bifurcation diagrams.
In the domain of smooth dynamical systems, the unfolding of the most common
codimension-two points is well known (see, e.g., \citep{Kuznetsov04}), and
this knowledge is exploited in continuation software for the automatic
switching among bifurcation branches at these points
(see, e.g., \citep{Dhooge03,Meijer09}).
The same cannot be said for nonsmooth systems, where,
though efficient numerical tools for bifurcation analysis are finally
starting to appear \citep{Dercole05c,Thota08b},
results are still mostly limited to codimension-one cases.
A reason for this shortcoming can be found in the fact that nonsmooth systems
exhibit, along with the standard bifurcations of smooth systems, a great
number of completely new bifurcations,
called \textit{discontinuity induced bifurcations},
that involve the interaction of the system's invariant sets with the
discontinuity boundaries. Since the characteristics of these bifurcations
depend critically on both the class of nonsmooth system and the geometry of
the involved boundaries, the number of possible scenarios is huge and,
at the moment, truly general results are scarce.
It goes without saying that codimension-two cases involving simultaneous
smooth and discontinuity induced bifurcations, named ``type II'' in \citep{Kowalczyk06}, are even more numerous, and less understood.

In this article we analyse type II bifurcations of periodic orbits (limit cycles), that is, bifurcations involving a periodic orbit (from now on called the bifurcating cycle) that collides with a discontinuity boundary while being at the same time non-hyperbolic.  Rather than aiming at a complete unfolding with reference to a particular class of nonsmooth systems, we concentrate on finding those geometric features that are common to all classes:
this is accomplished by abstracting our analysis from the nature of the involved boundary. As a consequence, our results are incomplete, because they only focus on the geometry of bifurcation curves around the codimension-two point; on the other hand, they apply more in general --- a feature that should be welcome in a field where peculiarity seems to be the rule.  

In particular, we show that three codimension-one bifurcation curves generically emanate from a type II point in a two-parameter plane.  One is the smooth bifurcation curve (fold, flip, or Neimark-Sacker), while the other two are the discontinuity induced bifurcations of the bifurcating cycle and of the secondary invariant set involved in the smooth bifurcation (the other cycle, the double-period cycle, or the torus, respectively).  Then we show that, depending on the bifurcation, one or both of these curves are tangent to the smooth bifurcation curve.  Indeed, in the flip and Neimark-Sacker cases, the bifurcating cycle departs from the image of the nonhyperbolic cycle, left frozen in state space, at a linear rate with respect to the bifurcation parameter, whereas the distance between the period-two cycle or the torus from such an image goes as the square root of the parameter perturbation from the bifurcation.  As a consequence, locally to the codimension-two point, the perturbation required by the secondary invariant set to collide with the discontinuity boundary is quadratic with respect to that required by the bifurcating cycle.  Similarly, in the fold case, the rate at which both cycles approach the image of the nonhyperbolic cycle is proportional to the square root of the parameter perturbation, so that the discontinuity induced bifurcation curves are both quadratically tangent to the fold curve.  These rather intuitive results have been observed in many examples, and proved for some specific classes of discontinuous systems (e.g., in \citep{Dankowicz05, Kowalczyk06, Nordmark06, Thota06, Zhao06, Simpson08b, Simpson09a}).  The aim of this paper is to provide formal support to the above geometric arguments and to prove their validity once and for all under very general conditions.

The ensuing exposition is set into the framework of
{\it grazing} bifurcations in continuous-time, where the discontinuity boundary
is smooth, locally to the point of contact with the bifurcating cycle, and
the contact occurs tangentially.
This allows us to keep the terminology as coherent as possible, especially
in the lack of a uniform terminology across all classes of nonsmooth systems.
Nonetheless, the reader will realise that our exposition is general and applies
to any discontinuity induced bifurcation involving a non-hyperbolic cycle in
continuous time or a non-hyperbolic fixed point in discrete time.
In fact, our analysis is based on the reduction of the nonsmooth flow to a
map which is defined and smooth on one side of a boundary, while we do not
describe the behaviour of the map on the other side. The rest of the analysis
is based on the obtained map, as if the problem was originally set in 
discrete time.
Thus, in practise, we do not make any assumption on the class of nonsmooth
systems and on the geometry of the discontinuity boundary.

We begin by stating the problem, introducing the basic notation, and
outlining the steps that we follow in the main proofs
(Sect.~\ref{sec:fa}); then we proceed with the detailed analysis of the three generic grazing bifurcations of non-hyperbolic cycles: the grazing-fold, the grazing-flip,
and the grazing-Neimark-Sacker (Sects.~\ref{sec:bfd}--\ref{sec:bns} and Appendices).
Once casted in discrete time, grazing bifurcations are more appropriately
called {\it border collisions}, and this is the name we use in this part
of the paper. Then we presents three specific applications
(Sect.~\ref{sec:ex}) and conclude with some future directions. 

\section{The framework of analysis}
\label{sec:fa}
We consider a nonsmooth autonomous flow
$x(t)=\Phi(x(0),t,\alpha)\in\R^{n+1}$
depending on parameters $\alpha\in\R^2$.
Namely, the right-hand side of the system's ODEs
\begin{equation}
\label{eq:ode}
\ds\dot{x}(t)=\left.\Frac{\partial}{\partial\tau}
\Phi(x(t),\tau,\alpha)\right|_{\tau=0}=
\Phi_t(x(t),0,\alpha)
\end{equation}
(here and in the following variables and parameters as subscripts denote
differentiation)
is generically smooth, but characterised by zero- or higher-order
discontinuities across some discontinuity boundaries $\mathcal{D}_i$,
defined as the zero set of suitable smooth functions $D_i(x,\alpha)$.
In particular, we can distinguish three types of discontinuity boundaries
(see Fig.~\ref{fig:cyc}):
\begin{figure}[t!]
\centerline{\includegraphics[width=1\textwidth]
{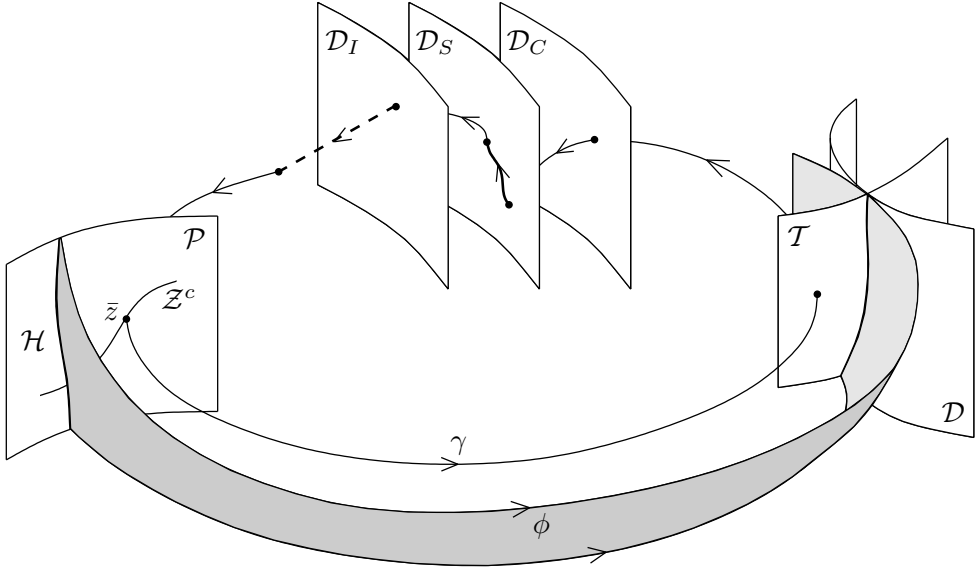}}
\caption{A generic (hyperbolic) limit cycle $\gamma$ of the nonsmooth flow $\Phi$.
For some $\alpha$ near $\alpha=0$, the cycle passes close to, but does not
touch, the discontinuity boundary $\mathcal{D}$, so that the resulting
Poincar\'e map on $\mathcal{P}$ is defined, locally to $\bar{z}$, only on
one side of the discontinuity boundary $\mathcal{H}$.
The boundary $\mathcal{H}$ divides $\mathcal{P}$ into two regions,
respectively composed of points $z$ from which the orbit of
$\Phi$ does and does not touch $\mathcal{D}$.}
\label{fig:cyc}
\end{figure}
boundaries across which the right-hand side of (\ref{eq:ode}) is
nonsmooth but continuous, so that orbits always cross the boundary
($\mathcal{D}_C$ in the figure);
boundaries across which the right-hand side of (\ref{eq:ode}) is
discontinuous, so that sliding motions are possible ($\mathcal{D}_S$);
boundaries where the right-hand side of (\ref{eq:ode}) is formally
characterised by impulsive components, which define an instantaneous state
transition (or jump) whenever orbits reach the boundary ($\mathcal{D}_I$).

Forward solutions of system (\ref{eq:ode}) are composed of smooth segments,
each corresponding to a smooth orbit terminating at a discontinuity boundary,
or to a sliding motion. Smooth segments are directly connected at crossing
and sliding boundaries, while they are connected through state jumps at
impacting boundaries.
Let $\gamma$ be a periodic orbit of system (\ref{eq:ode}).
In Fig.~\ref{fig:cyc}, $\gamma$ is composed of four segments, three smooth
(solid) orbits and one sliding motion (thick orbit),
and is characterised by a single state jump (thick dashed connection).

Suppose that, when $\alpha=0$, the cycle $\gamma$ grazes (touches tangentially) a
discontinuity boundary $\mathcal{D}$, and no other degeneracies occur on
$\mathcal{D}_C$, $\mathcal{D}_I$, and $\mathcal{D}_S$.  At the same time, suppose that $\gamma$  is non-hyperbolic at $\alpha=0$ (more precisely, the multipliers are not defined at $\alpha=0$, but the smooth bifurcation curve is a path to $\alpha=0$ on which one real or two complex conjugate simple multipliers lie on the unit circle).
Introduce a Poincar\'e section $\mathcal{P}$ along one of the segments of
$\gamma$, say, e.g., the segment touching $\mathcal{D}$ so that the
flow reaches $\mathcal{D}$ after $\mathcal{P}$ for $\alpha=0$.
Also introduce a coordinate $z\in\R^n$ on $\mathcal{P}$ such that the
intersection $\bar{z}$ of $\gamma$ with $\mathcal{P}$ lies at $z=0$
for $\alpha=0$.
Then, locally to $(z,\alpha)=(0,0)$, the flow $\Phi$ induces a Poincar\'e map
\begin{equation}
\label{eq:pm}
z\mapsto F(z,\alpha)
\end{equation}
(note that the map may not be invertible, e.g., in the presence of sliding motions).  Since we do not discuss the type of boundary $\mathcal{D}$, we limit the definition of $F$ to the values of $(z,\alpha)$ in a neighbourhood of $(0,0)$ for which the orbit originating at $z$ does not touch $\mathcal{D}$.
This introduces an $(n-1)$-dimensional discontinuity boundary $\mathcal{H}$ on the Poincar\'e section $\mathcal{P}$ such that $F$ is defined and smooth on one side of $\mathcal{H}$. In particular, let
$$
\mathcal{D}=\{x:D(x,\alpha)=0\},\quad \mathcal{H}=\{z:H(z,\alpha)=0\},
$$
and assume, without loss of generality, that the flow $\Phi$ touches
$\mathcal{D}$ tangentially while locally remaining on the side $D(x,\alpha)<0$,
and that $F(z,\alpha)$ is defined for $H(z,\alpha)<0$.
Then, the function $H$ can be constructed as follows
(see again Fig.~\ref{fig:cyc}).
Define the $n$-dimensional smooth manifold $\mathcal{T}$ of the points where
the flow is tangent to the level sets of function $D$:
$$
\mathcal{T}=\{x:T(x,\alpha):=
\left\langle\Phi_t(x,0,\alpha),D_x(x,\alpha)\right\rangle=0\}
$$
(vector $D_x(x,\alpha)\in\R^{n+1}$ is orthogonal to the level sets of
$\mathcal{D}$ at $(x,\alpha)$ and $\langle\cdot,\cdot\rangle$ is the standard
scalar product in $\R^{n+1}$).
As shown in Fig.~\ref{fig:cyc},
the $(n-1)$-dimensional intersection between $\mathcal{D}$ and $\mathcal{T}$
is transformed, backward in time by the flow, into the discontinuity boundary
$\mathcal{H}$.
Thus, $H(z,\alpha)$ can be defined as the value $D(x,\alpha)$ at
the point $x$ at which the flow first reaches $\mathcal{T}$ (forward in time)
from the initial condition corresponding to $z$ on $\mathcal{P}$.

We can now abandon the continuous-time framework, and focus on
map (\ref{eq:pm}).
For some $\alpha$ in a neighbourhood of $\alpha=0$, the map is characterised
by a fixed point $\bar{z}$, with $H(\bar{z},\alpha)<0$ and, for $\alpha=0$,
the fixed point is non-hyperbolic and lies at the origin $z=0$ and on the
discontinuity boundary $\mathcal{H}$.
We investigate the bifurcation curves rooted at $\alpha=0$ in the
parameter plane $(\alpha_1,\alpha_2)$,
by considering separately the three generic cases, namely
(I) fold (one simple eigenvalue equal to $1$, Sect.~\ref{sec:bfd}),
(II) flip (one simple eigenvalue equal to $-1$, Sect.~\ref{sec:bfp}), and
(III) Neimark-Sacker
(two simple complex conjugate eigenvalues on the unit circle,
Sect.~\ref{sec:bns}).

In each case, we proceed as follows. Locally to $(z,\alpha)=(0,0)$, we
consider the restriction of map (\ref{eq:pm}) to a parameter-dependent
centre manifold $\mathcal{Z}^c$.
Let $u\in\R^{n_c}$ represent coordinates on $\mathcal{Z}^c$, $n_c=1$ in the
fold and flip cases, $n_c=2$ in the Neimark-Sacker case, with $u=u(z,\alpha)$
for each $z\in\mathcal{Z}^c$ and $\alpha$ in a neighbourhood of
$(z,\alpha)=(0,0)$, $u(0,0)=0$, and let $z=z(u,\alpha)$ denote the inverse
transformation. Restricted to the centre manifold, map (\ref{eq:pm}) reads
\begin{equation}
\label{eq:cm}
u\mapsto f(u,\alpha):=u(F(z(u,\alpha),\alpha),\alpha).
\end{equation}
and the discontinuity boundary $\mathcal{H}$ is given by the zero-set of the function
\begin{equation}
\label{eq:h_func}
h(u,\alpha):=H(z(u,\alpha),\alpha).
\end{equation}
We assume that the three following conditions hold:
\begin{itemize}
 \item[(i)] Map (\ref{eq:cm}) satisfies, at $\alpha=0$, all genericity conditions of the corresponding smooth bifurcation (see, e.g., \citep{Kuznetsov04}).
 \item[(ii)] At $\alpha=0$, the centre manifold $\mathcal{Z}^c$ transversely intersects the
discontinuity boundary $\mathcal{H}$ at $z=0$ (by continuity the transversality persists near $(z,\alpha)=(0,0)$, see Fig.~\ref{fig:cyc}).  Under this condition, the dynamics of map (\ref{eq:pm}) near $(z,\alpha)=(0,0)$ is captured by that on the centre manifold.  In the coordinate $u$ along the centre manifold the condition becomes $h_u(0,0)\neq 0$.
 \item[(iii)] Changing $\alpha$ along the smooth bifurcation curve, the non-hyperbolic fixed point crosses the discontinuity boundary transversely.  This condition ensures that the smooth bifurcation curve intersects the border collision curves in a generic way.
\end{itemize}
As a first step, we reduce map (\ref{eq:cm}) to a normal form (NF)
(the fold, flip, and Neimark-Sacker normal forms)
through a locally invertible change of variable and parameter, say, 
$v=v(u,\alpha)$, $\beta=\beta(\alpha)$, where $v(0,0)=0$, $\beta(0)=0$, and
$u=u(v,\beta)$, $\alpha=\alpha(\beta)$ denote the inverse transformation.
Then, second step, we find the expression of the discontinuity boundary (\ref{eq:h_func}) 
in the new variables and parameters, i.e.,
\begin{equation}
\label{eq:Sigma}
\{v:h^{\mathrm{NF}}(v,\beta):=h(u(v,\beta),\alpha(\beta))=0\}. 
\end{equation}
Finally, third step, we analyse the interaction of the normal form map
$$
v\mapsto f^{\mathrm{NF}}(v,\beta):=
v(f(u(v,\beta),\alpha(\beta)),\alpha(\beta))
$$
with the discontinuity boundary (\ref{eq:Sigma}), and we find local asymptotics for the bifurcation curves emanating from $\alpha=0$ in terms of $(\alpha_1,\alpha_2)$-expansions.

The details of the normal form reduction are reported in appendices A.1, B.1, and C.1, while the technicalities on step two are reported in Appendices A.2, B.2, and C.2.  The specific analytical form taken by condition (iii) in the fold, flip, and NS cases is respectively derived in Appendices A.3, B.3, and C.3 in terms of both the original coordinates $z$ and in the coordinates $u$ in the centre manifold. Finally, some details on step three for the Neimark-Sacker case are relegated to Appendix C.4. For simplicity of notation, in the following the 0 superscript stands for evaluation at $(u,\alpha)=(0,0)$ or $(v,\beta)=(0,0)$.

\section{Case I: Border-fold bifurcation}
\label{sec:bfd}
Let the dynamics in the centre manifold $\mathcal{Z}^c$ be described by the
one-dimensional system
\begin{equation}
\label{eq:bfd_cm}
u \mapsto f(u,\alpha), \quad u\in\R^1,
\end{equation}
with $f^0=0$ (fixed point condition) and $f_u^0=1$ (fold condition).  Under condition (i), map (\ref{eq:bfd_cm}) can be reduced to normal form (first step, see Appendix A.1) with invertible changes of variable and parameter  $v=v(u,\alpha)$, $\beta=\beta(\alpha)$, becoming
\begin{equation}
\label{eq:bfd_nf}
v\mapsto\beta_1+v+sv^2+O(v^3),
\end{equation}
where $s=\mathrm{sign}(f^0_{uu})$.  In these variables, the fold curve has equation $\beta_1=0$ in the plane $(\beta_1,\beta_2)$, and the corresponding non-hyperbolic fixed point is located at $v=0$.

We now turn our attention to the discontinuity boundary (\ref{eq:Sigma}) (second step, see Appendix A.2).
Condition (ii), ensuring transversal intersection of the centre manifold $\mathcal{Z}^c$ and the discontinuity boundary $\mathcal{H}$, implies local existence and uniqueness of a smooth function
$$
\sigma(\beta)=\sigma_{\beta_1}^0\beta_1+\sigma_{\beta_2}^0\beta_2+
O(\|\beta\|^2),
$$
such that the intersection of $\mathcal{H}$
with $\mathcal{Z}^c$ is located at $v=\sigma(\beta)$.  Then by condition (iii) (see Appendix A.3 for the analytical expression) we know that moving along the fold curve, that is, along the $\beta_2$-axis, the fixed point at $v=0$ crosses $\mathcal{H}$ at $\beta_2=0$.  As a consequence, we have $\sigma_{\beta_2}^0\neq0$.

We are now ready to find the equation of the border collisions in the plane $(\beta_1,\beta_2)$ (third step). The two fixed points of the normal form map (\ref{eq:bfd_nf}) are located at $\bar{v}^{\pm}(\beta)=\pm\sqrt{-s\beta_1}+O(\|\beta\|^2)$ ($\bar{v}^-$ being stable and $\bar{v}^+$ unstable for $s=1$, and viceversa for $s=-1$), and lie on the
discontinuity boundary (\ref{eq:Sigma}) along the curves 
\begin{equation}
\label{eq:bfd_1}
\pm\sqrt{-s\beta_1}=\sigma_{\beta_1}^0\beta_1+\sigma_{\beta_2}^0\beta_2+
O(\|\beta\|^2).
\end{equation}
Since $\sigma_{\beta_2}^0\neq 0$,
equation (\ref{eq:bfd_1}) for small $\|\beta\|$ becomes
\begin{equation}
\label{eq:bfd_2}
\pm\sqrt{-s\beta_1}\simeq \sigma_{\beta_2}^0\beta_2,
\end{equation}
and gives the asymptotics, locally to $\beta=0$, of the two
border-collision bifurcation curves involving the fixed points $\bar{v}^{\pm}$.  The invertible parameter change $\beta=\beta(\alpha)$ easily provides the asymptotics in the original $\alpha$ parameters.

\begin{figure}[t!]
\centerline{\includegraphics[width=1\textwidth]
{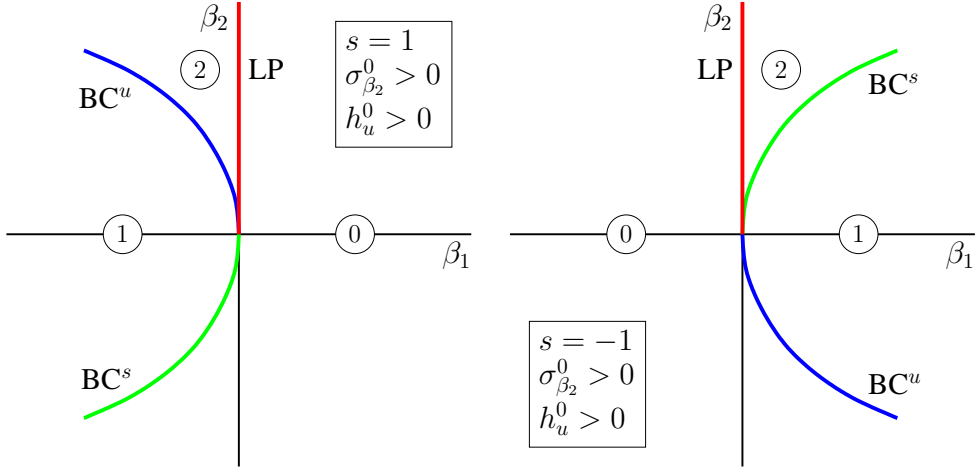}}
\caption{Border-fold bifurcation.
Bifurcation curves: LP, fold (limit point, red);
BC$^{s}$, border collision of the stable fixed point ($\bar{v}^{-}$, left; $\bar{v}^{+}$, right) of map (\ref{eq:bfd_nf}) (green), BC$^{u}$, border collision of the unstable fixed point ($\bar{v}^{+}$, left; $\bar{v}^{-}$, right) of map (\ref{eq:bfd_nf}) (blue).
Region labels: $0$, no fixed point in $V^{-}(\beta):=\{v:h^{\mathrm{NF}}(v,\beta)<0\}$;
$1$, $\bar{v}^{-}$ is the only fixed point in $V^{-}(\beta)$;
$2$, both fixed points $\bar{v}^{\pm}$ lie in $V^{-}(\beta)$.}
\label{fig:bfd}
\end{figure}

Depending upon the sign of $s$ in the normal form map (\ref{eq:bfd_nf}), of $\sigma_{\beta_2}^0$ in (\ref{eq:bfd_2}), and of $h_u^0$ in (ii), there are eight generic cases, two of which are reported in Fig.~\ref{fig:bfd}.
The other six can be reduced to these two by suitable parameter changes.  In fact, the four cases with $\sigma_{\beta_2}^0<0$ are symmetric with respect to the $\beta_1$-axis to the corresponding cases with $\sigma_{\beta_2}^0>0$, while the four cases with $h_u^0<0$ can be reduced to cases with $h_u^0>0$ by changing the sign of $s$ and rotating the figure.  Note that only half of the $\beta_2$-axis can be said to belong to the fold curve (LP), since along the other half the two fixed points $\bar{v}^{\pm}$
collide at $v=0$ on the undescribed side of the discontinuity boundary
(\ref{eq:Sigma}), i.e., $h^{\mathrm{NF}}(0,\beta)>0$.

\section{Case II: Border-flip bifurcation}
\label{sec:bfp}
Let the dynamics in the centre manifold $\mathcal{Z}^c$ be described by the
one-dimensional system
\begin{equation}
\label{eq:bfp_cm}
u \mapsto f(u,\alpha), \quad u\in\R^1,
\end{equation}
with $f^0=0$ (fixed point condition) and $f_u^0=-1$ (flip condition).  Through a parameter-dependent translation, we can ensure that $f(0,\alpha)=0$, i.e., that $u=0$ is a fixed point for all $\alpha$ in a neighbourhood of $\alpha=0$.   Under condition (i), map (\ref{eq:bfp_cm}) can be reduced to normal form (first step, see Appendix B.1) with invertible changes of variable and parameter  $v=v(u,\alpha)$, $\beta=\beta(\alpha)$, becoming
\begin{equation}
\label{eq:bfp_nf}
v\mapsto-(1+\beta_1)v + sv^3+O(v^4),
\end{equation}
with $s=\mathrm{sign}((1/4)(f^0_{uu})^2+(1/6)f^0_{uuu})$.  In these variables, the flip curve has equation $\beta_1=0$ in the plane $(\beta_1,\beta_2)$, and the corresponding non-hyperbolic fixed point is located at $v=0$.  Moreover, parameters can be chosen so that the border collision of the fixed point in the origin has equation $\beta_2=0$.

We now turn our attention to the discontinuity boundary (\ref{eq:Sigma}) (second step, see Appendix B.2).
Condition (ii), ensuring transversal intersection of the centre manifold $\mathcal{Z}^c$ and the discontinuity boundary $\mathcal{H}$, implies local existence and uniqueness of a smooth function
$$
\sigma(\beta)=\sigma_{\beta_1}^0\beta_1+\sigma_{\beta_2}^0\beta_2+
O(\|\beta\|^2),
$$
such that the intersection of $\mathcal{H}$ with $\mathcal{Z}^c$ is located at $v=\sigma(\beta)$.  Moreover, thanks to the parameter choice in (\ref{eq:bfp_nf}), $\sigma_{\beta_1}^0=0$, since the fixed point $v=0$ lies on $\mathcal{H}$ when $\beta_2=0$.  Then by condition (iii) (see Appendix B.3 for the analytical expression) we know that moving along the flip curve, that is, along the $\beta_2$-axis, the fixed point at $v=0$ crosses $\mathcal{H}$ at $\beta_2=0$.  As a consequence, we have $\sigma_{\beta_2}^0\neq0$.

We are now ready to find the equation of the border collisions in the plane $(\beta_1,\beta_2)$ (third step).  Near $(v,\beta_1)=(0,0)$ the normal form map (\ref{eq:bfp_nf}) iterated twice
has one fixed point in $v=0$ (which is also a fixed point of map (\ref{eq:bfp_nf}))
and two others in $\bar{v}^{\pm}(\beta)=\pm\sqrt{s\beta_1}+O(\|\beta\|^2)$
(period-two cycle). In particular, $\bar{v}^{\pm}$ lie on discontinuity boundary (\ref{eq:Sigma}) along the curves
\begin{equation}
\label{eq:bfp_1}
\pm\sqrt{s\beta_1}=\sigma_{\beta_2}^0\beta_2+
O(\|\beta\|^2).
\end{equation}
Since $\sigma_{\beta_2}^0\neq 0$,
equation (\ref{eq:bfp_1}) for small $\|\beta\|$ becomes
\begin{equation}
\label{eq:bfp_2}
\pm\sqrt{s\beta_1}\simeq \sigma_{\beta_2}^0\beta_2,
\end{equation}
and gives the asymptotics, locally to $\beta=0$, of the
border-collision bifurcation curves involving the two points $\bar{v}^{\pm}$
of the period-two cycle.  The invertible parameter change $\beta=\beta(\alpha)$ provides the asymptotics in the original $\alpha$ parameters.

\begin{figure}[t!]
\centerline{\includegraphics[width=1\textwidth]
{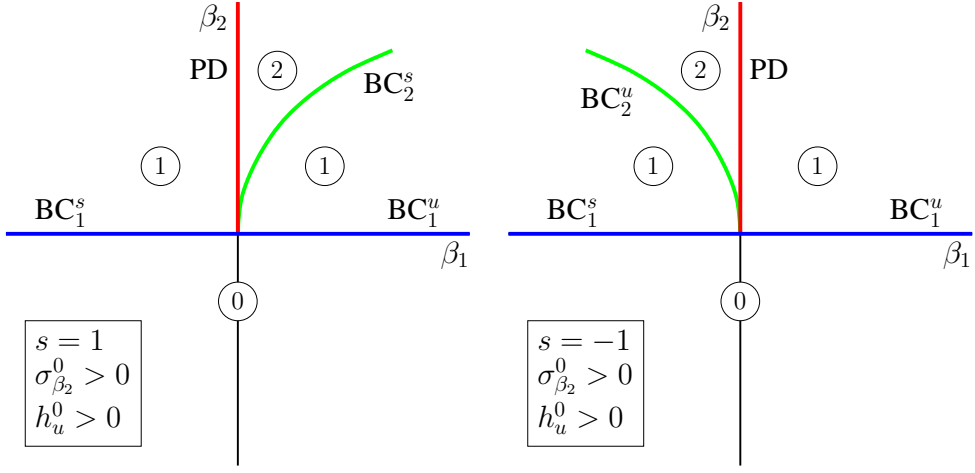}}
\caption{Border-flip bifurcation.
Bifurcation curves: PD, flip (period doubling, red);
BC$_1^{s,u}$, border collision of the fixed point $v=0$ (stable and unstable branches, blue);
BC$_2^{s,u}$, border collision of the stable or unstable period-two cycle.
Region labels: $0$, no fixed point or period-two cycle in $V^{-}(\beta):=\{v:h^{\mathrm{NF}}(v,\beta)<0\}$;
$1$, $v=0$ is a fixed point in $V^{-}(\beta)$ and there is no period-two cycle,
or it does not lie entirely in $V^{-}(\beta)$;
$2$, the fixed point $v=0$ coexists in $V^{-}(\beta)$ with the period-two cycle.}
\label{fig:bfp}
\end{figure}

Depending upon the sign of $s$ in the normal form map
(\ref{eq:bfp_nf}), of $\sigma_{\beta_2}^0$ in
(\ref{eq:bfp_2}), and of $h_u^0$ in (ii),
there are eight generic cases. However, again, only two cases are relevant
(see Fig.~\ref{fig:bfp}),
because all others can be reduced to these two by
suitable parameter changes.  Here, both the four cases with $\sigma_{\beta_2}^0<0$ and those with $h_u^0<0$, are symmetric with respect to the $\beta_1$-axis to the corresponding cases with $\sigma_{\beta_2}^0>0$ or $h_u^0>0$.  Also note that only half of the $\beta_2$-axis can be said to belong to the flip curve (PD), since along the other half the fixed point $v=0$ lies on the
undescribed side of the discontinuity boundary (\ref{eq:Sigma}), i.e., $h^{\mathrm{NF}}(0,\beta)>0$.
Similarly, only one of the two branches in (\ref{eq:bfp_2}) constitutes the border-collision curve involving the period-two cycle (stable, BC$_2^s$; unstable, BC$_2^u$), since along the other branch $h^{\mathrm{NF}}(\bar{v}^\pm,\beta)\geq0$.

\section{Case III: Border-Neimark-Sacker bifurcation}
\label{sec:bns}
Let the dynamics in the centre manifold $\mathcal{Z}^c$ be described by the
two-dimensional system
\begin{equation}
\label{eq:bns_cm}
u \mapsto f(u,\alpha), \quad u\in\R^2,
\end{equation}
with $f^0=0$ (fixed point condition) and with eigenvalues
$\lambda^0$ and $\bar{\lambda}^0$ 
(the overbar stands for complex conjugation)
of the $2\times 2$ Jacobian $f_u^0$ given by
$$
\lambda(\alpha)=(1+g(\alpha))e^{i\theta(\alpha)},
$$
with $g^0=0$ (Neimark-Sacker, NS, condition).  As in the flip case, assume that $f(0,\alpha)=0$ for all $\alpha$ in a neighbourhood of $\alpha=0$.  Under condition (i), map (\ref{eq:bns_cm}) can be reduced to normal form in polar coordinates (first step, see Appendix C.1) with invertible changes of variable and parameter  $\rho=\rho(u,\alpha), \varphi=\varphi(u,\alpha)$, $\beta=\beta(\alpha)$, becoming
\begin{subequations}
\label{eq:bns_nf}
\begin{eqnarray}
\label{eq:bns_nf_r}
\rho & \mapsto & \rho(1+\beta_1+a(\beta)\rho^2)+\rho^4 R(\rho,\varphi,\beta),\\
\label{eq:bns_nf_p}
\varphi & \mapsto & \varphi+\theta(\alpha(\beta))+
\rho^2 Q(\rho,\varphi,\beta),
\end{eqnarray}
\end{subequations}
where $a^0\neq 0 $. In these variables, the NS curve has equation $\beta_1=0$ in the plane $(\beta_1,\beta_2)$, and the corresponding non-hyperbolic fixed point is located at $v=0$ (with $v_1=\mathrm{Re}(\rho e^{i\varphi})$ and $v_2=\mathrm{Im}(\rho e^{i\varphi})$).  Moreover, parameters can be chosen so that the border collision of the fixed point in the origin has equation $\beta_2=0$.

We now turn our attention to the discontinuity boundary (\ref{eq:Sigma}) (second step, see Appendix C.2).  Condition (ii), ensuring transversal intersection of the centre manifold $\mathcal{Z}^c$ and the discontinuity boundary $\mathcal{H}$, implies local existence and uniqueness of a smooth function
$$
\sigma(\beta)=\sigma_{\beta_1}^0\beta_1+\sigma_{\beta_2}^0\beta_2+
O(\|\beta\|^2),
$$
measuring the distance between the origin and the boundary, with positive/negative values if $h^{\mathrm{NF}}(0,\beta)$ is negative/positive, in order to make $\sigma(\beta)$ differentiable at $\beta=0$.  Moreover, thanks to the parameter choice in (\ref{eq:bns_nf}), $\sigma_{\beta_1}^0=0$, since the fixed point $v=0$ lies on $\mathcal{H}$ when $\beta_2=0$. Then by condition (iii) (see Appendix C.3 for the analytical expression) we know that moving along the NS curve, that is, along the $\beta_2$-axis, the fixed point at $v=0$ crosses $\mathcal{H}$ transversely at $\beta_2=0$.  As a consequence, we have $\sigma_{\beta_2}^0\neq0$.

We are now ready to find the equation of the border collisions in the plane $(\beta_1,\beta_2)$ (third step). Near $\beta=0$, the normal form map (\ref{eq:bns_nf}) has an fixed point in
$\rho=0$ and a closed invariant curve that is contained in the annular region
\begin{equation}
\label{eq:bns_ar}
\left\{(\rho,\varphi):
\sqrt{-\Frac{\beta_1}{a(\beta)}}(1-\beta_1^{\gamma-1/2})\le\rho\le
\sqrt{-\Frac{\beta_1}{a(\beta)}}(1+\beta_1^{\gamma-1/2}),\,
\varphi\in[0,\,2\pi]\right\},\quad\Frac{1}{2}<\gamma<1
\end{equation}
(see Appendix C.4).  The two circles delimiting the annular region (\ref{eq:bns_ar})
touch the discontinuity boundary along the curves
\begin{equation}
\label{eq:bns_1}
\sqrt{-\Frac{\beta_1}{a(\beta)}}(1\pm\beta_1^{\gamma-1/2})=
\sigma_{\beta_2}^0\beta_2+O(\|\beta\|^2).
\end{equation}
Since $\sigma_{\beta_2}^0\neq 0$,
equation (\ref{eq:bns_1}) for small $\|\beta\|$ becomes
\begin{equation}
\label{eq:bns_2}
\sqrt{-\Frac{\beta_1}{a^0}} \simeq \sigma_{\beta_2}^0\beta_2,
\end{equation}
and gives a unique asymptotic, locally to $\beta=0$, for the grazing
bifurcation curves of both circles. The same asymptotic therefore holds
for the grazing bifurcation involving the invariant curve (the uniqueness of the bifurcation curve is granted by the elliptical shape of the invariant curve near $\beta=0$).
Again, the invertible parameter change $\beta=\beta(\alpha)$ provides the asymptotics in the original $\alpha$ parameters.

\begin{figure}[t!]
\centerline{\includegraphics[width=1\textwidth]
{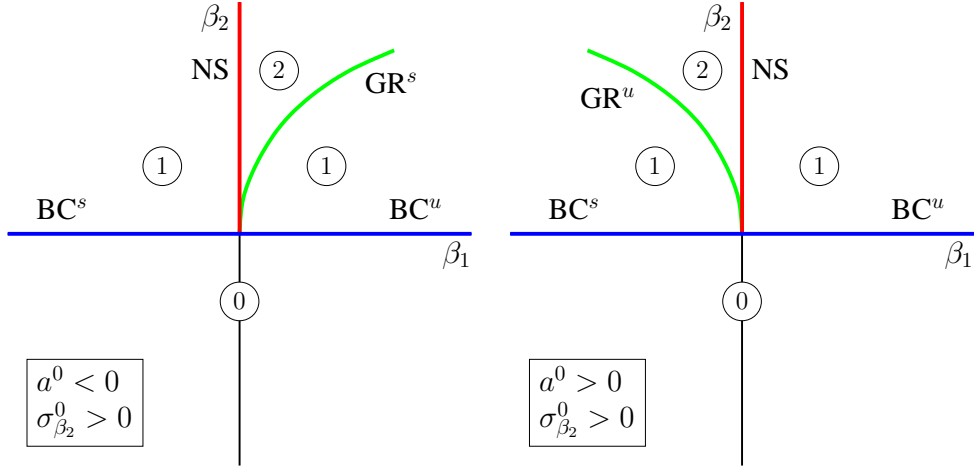}}
\caption{Border-NS bifurcation.
Bifurcation curves: NS, Neimark-Sacker (red);
BC$^{s,u}$, border collision of the fixed point $v=0$ (stable and unstable branches, blue);
GR$^{s,u}$, grazing of the stable or unstable torus (green).
Region labels: $0$, no fixed point or invariant curve in
$V^{-}(\beta):=\{v:h^{\mathrm{NF}}(v,\beta)<0\}$;
$1$, $v=0$ is a fixed point in $V^{-}(\beta)$ and there is no invariant curve,
or it does not lie entirely in $V^{-}(\beta)$;
$2$, both the fixed point $v=0$ and the invariant curve lie in $V^{-}(\beta)$.}
\label{fig:bns}
\end{figure}

Depending upon the sign of $a^0$ in the normal form map (\ref{eq:bns_nf}) and of $\sigma_{\beta_2}^0$ in (\ref{eq:bns_2}), there are four generic cases. However, again, only two cases are relevant (see Fig.~\ref{fig:bns}),
because those with $\sigma_{\beta_2}^0<0$ are symmetric with respect to the $\beta_1$-axis to the cases with $\sigma_{\beta_2}^0>0$.  Also note that only half of the $\beta_2$-axis can be said to belong to the
NS curve, since along the other half the fixed point $v=0$ lies on the
undescribed side of the discontinuity boundary (\ref{eq:Sigma}), i.e.,
$h^{\mathrm{NF}}(0,\beta)>0$.  Similarly, only half of the parabola in (\ref{eq:bns_2}) constitutes the grazing bifurcation curve involving the invariant curve (stable, GR$^s$; unstable, GR$^u$), since along the other half the invariant curve is composed of points $v$ with $h^{\mathrm{NF}}(v,\beta)\geq0$.

\section{Examples}
\label{sec:ex}
We now present three specific examples, one for each of the three
codimension-two bifurcations analysed in the previous sections.
The three examples deal with different classes of nonsmooth systems
(an impacting, a hybrid, and a piecewise smooth system) and describe interesting applications in different fields of science and engineering
(ecology, social sciences, and mechanics).

\subsection*{An impacting model of forest fires}
\label{sec:ex1}
For an example of border-fold bifurcation, we consider the
forest fire impacting model presented in \citep{Dercole05a,Maggi06}.
The model describes the vegetational growth with the following two
(smooth) ODEs:
\begin{eqnarray*}
\dot{B} & = & \ds r_B B \left(1 - \Frac{B}{K_B}\right) - \alpha B T,\\
\dot{T} & = & \ds r_T T \left(1 - \Frac{T}{K_T}\right),
\end{eqnarray*}
one for the surface layer (bush, $B$) and one for the upper layer (trees, $T$).
Fire episodes are represented by instantaneous events (impacts), that occur
when the biomasses $(B,T)$ of the two layers reach one of three
specified impacting boundaries:
a bush ignition threshold $\rho_BK_B$ triggering bush-only fires that map the bush biomass to $\lambda_B\rho_BK_B$, $0<\lambda_B,\rho_B<1$;
a tree ignition threshold $\rho_TK_T$ triggering trees-only fires that map the trees biomass to $\lambda_T\rho_TK_T$, $0<\lambda_T,\rho_T<1$;
and the segment connecting points $(\sigma_BK_B,\rho_TK_T)$ and
$(\rho_BK_B,\sigma_TK_T)$, $0<\sigma_B<\rho_B$, $0<\sigma_T<\rho_T$,
triggering mixed fires with post-fire conditions suitably assigned as a
function of pre-fire conditions (see \citep{Maggi06} for more details).

\begin{figure}[t!]
\centerline{\includegraphics[scale=0.9]
{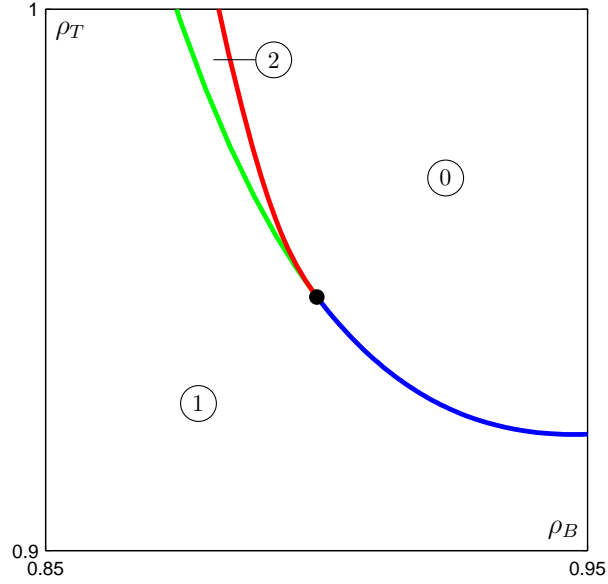}}
\caption{Example of border-fold bifurcation.
Bifurcation curves: fold (red);
border collision of the period-one stable cycle (blue);
border collision of the period-one unstable cycle (green).
Region labels as in Fig.~\ref{fig:bfd}.}
\label{fig:ex1}
\end{figure}

For the parameter setting $r_1=0.375$, $r_2=0.0625$, $\alpha=0.43$,
$K_B=K_T=1$, $\rho_B=0.85$, $\rho_T=0.93$, $\lambda_B=0.03$, $\lambda_T=0.01$,
$\sigma_B=0.61$, $\sigma_T=0.3$ (corresponding to Mediterranean forests),
the system is characterised by a globally stable period-one cycle composed
of a growth orbit and a mixed fire.
Numerical continuation (by means of {\sc Auto07p} \citep{Auto07p}) of the
cycle in the parameter plane $(\rho_B,\rho_T)$ identifies two (codimension-one)
bifurcations: a fold (red curve in Fig.~\ref{fig:ex1})
and a grazing of the
growth orbit with the bush ignition threshold (blue curve).
The two curves merge together at the border-fold bifurcation (black) point
and, as predicted by the analysis carried out in Sect.~\ref{sec:bfd}, the
grazing bifurcation of the unstable cycle involved in the fold (green curve)
emanates tangentially to the fold curve from the codimension-two bifurcation point.

\subsection*{A hybrid model of two-party democracies}
\label{sec:ex2}
For an example of border-flip bifurcation, we consider the hybrid model
presented in \citep{Colombo08a} for describing the dynamics of
two-party democracies.  The model describes the evolution of the size of two lobbies (of sizes $L_D$ and $L_R$), one associated to each party (parties $D$ and $R$, respectively), and assumes that the individuals belonging to the lobby of the party at the government erode the welfare ($W$) at a rate proportional to the size of the lobby; a lobby can grow only as long as its party is at the government, and decays otherwise; a small fraction of the lobbyists not at the government defect and switch to the other lobby; elections are held once every $T$ years, and people vote for the party that has the less damaging lobby at the time of the elections.  Altogether, the dynamics is captured by two sets of ODEs,
namely
\begin{eqnarray*}
\dot{W} & = & r(1-W-a_DL_D)W,\\
\dot{L}_D & = & (e_Da_DW-d_D)L_D + k_RL_R,\\
\dot{L}_R & = & (-d_R-k_R)L_R,
\end{eqnarray*}
when the $D$-party is at the government, and
\begin{eqnarray*}
\dot{W} & = &  r(1-W-a_RL_R)W,\\
\dot{L}_D & = & (-d_D-k_D)L_D,\\
\dot{L}_R & = & (e_Ra_RW-d_R)L_R + k_DL_D,
\end{eqnarray*}
when the $R$-party is at the government.  Here, $r$ is the intrinsic growth rate of the welfare, $a$ represents the aggressiveness of a lobby, $e$ is the recruitment coefficient of a lobby, and $d$ and $k$ are respectively the rate at which individuals abandon the lobbies or defect. In the region of the state space where $a_DL_D < a_RL_R$ ($a_DL_D > a_RL_R$) the $D$-lobby ($R$-lobby) is less damaging and thus wins the elections. The condition $a_DL_D = a_RL_R$ therefore defines the discontinuity boundary (see \citep{Colombo08a} for more details).

\begin{figure}[t!]
\centerline{\includegraphics[scale=0.9]
{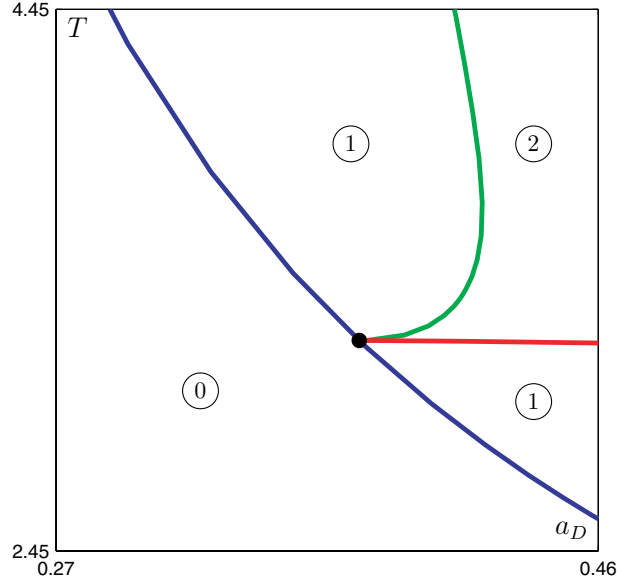}}
\caption{Example of border-flip bifurcation.
Bifurcation curves: flip (red);
border collision of the period-one cycle (blue);
border collision of the period-two cycle (green).
Region labels as in Fig.~\ref{fig:bfp}.
}
\label{fig:ex2}
\end{figure}

In the $(a_D,T)$ plane, with parameters $a_R=1$, $r=0.2$, $e_D=e_R=6$, $d_D=d_R=1.8$, $k_D=k_R=0.06$, the system has a very complex bifurcation diagram (see for example Fig.~1 in \citep{Colombo08a}).  In particular, near $a_D=0.38$, $T=3.2$, a flip (red curve in Fig.~\ref{fig:ex2})
and a border collision (blue curve) of a period-$2T$ cycle meet at the border-flip (black) point and, as predicted by the analysis carried out in Sect.~4, a border collision of the period-4T cycle (green curve) emanates from the codimension-two point tangentially to the flip curve.

\subsection*{A piecewise smooth model of railway wheelset dynamics}
\label{sec:ex3}
For an example of border-NS bifurcation, we consider a two degrees of freedom piecewise smooth model of a suspended railway wheelset with dry friction dampers, subject to a sinusoidal disturbance representing the deformations of the track.  The model is based on that presented in \citep{True03, Knudsen92}, where the track deformation was not taken into account, and its analysis will be published elsewhere.  Since a detailed explanation of the equations and parameters  goes beyond the scope of this paper, here we only report the equations and describe a few key parameters (see \citep{Knudsen92} and \citep{True03} for the details).   The model consists of the following piecewise smooth equations:
$$
\begin{array}{rcl}
 \dot{x}_1 & = & \tilde{x}_2,\\
 \dot{x}_2 & = & \Frac{1}{m}(-2F_x-2K_s \tilde{x}_1-	\mbox{sign}(x_2)\mu),\\
 \dot{x}_3 & = & x_4,\\
 \dot{x}_4 & = & \Frac{1}{I}(-2AF_y),
\end{array}
$$
where
$$
\begin{array}{c}
\tilde{x}_1=x_1+a\sin(\omega t),\quad \tilde{x}_2=x_2+a\omega\cos(\omega t),\\ \mu=(\mu_d(1-\sech(\alpha \tilde{x}_2)) + \mu_s \sech(\alpha \tilde{x}_2)),\\
F_x=\Frac{\xi_x F_r}{\Psi \xi_r}, \quad F_y=\Frac{\xi_y F_r}{\Phi \xi_r},\quad F_r=
\left\{\begin{array}{l}
\xi_r C\left(1-\Frac{C \xi_r}{3 \mu_t} + \Frac{C^2 \xi^2_r}{27 \mu^2_t}\right)\mbox{ if $C\xi_r<3\mu_t$},\\
\mu_t \mbox{ otherwise,}\\
\end{array}\right.\\
\xi_x=\Frac{\tilde{x}_2}{V}-x_3, \quad \xi_y=\Frac{Ax_4}{V}+\Frac{\lambda \tilde{x}_1}{r_0}, \quad \xi_r=\sqrt{\left(\Frac{\xi_x}{\Psi}\right)^2 + \left(\Frac{\xi_y}{\Phi}\right)^2}.
\end{array}
$$
Here $\omega=2\pi V/l$, $a$ and $l$ are the amplitude and wavelength of the sinusoidal disturbance, $V$ is the speed of the wheelset, and $\lambda$ measures the conicity of the wheels. The system's state space is therefore partitioned in four regions, depending on the signs of $x_2$ and of $C\xi_r-3\mu_t$, so that $x_2=0$ and $C\xi_r=3\mu_t$ define two discontinuity boundaries.  

\begin{figure}[t!]
\centerline{\includegraphics[scale=0.9]
{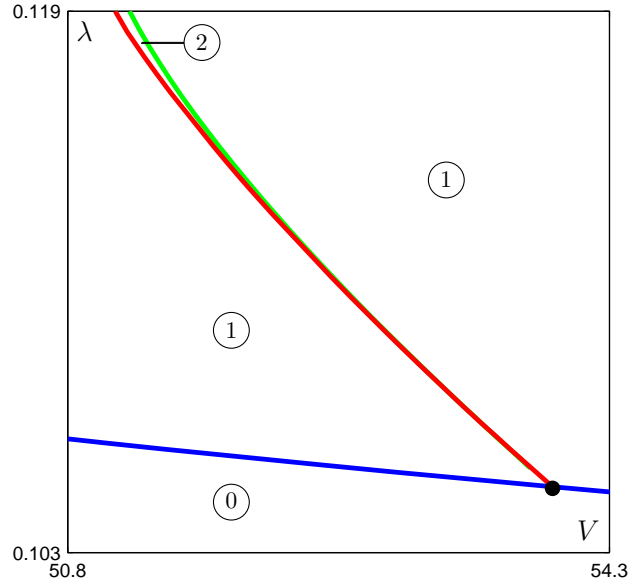}}
\caption{Example of border-NS bifurcation.
Bifurcation curves: Neimark-Sacker (red);
border collision of the period-one cycle (blue);
border collision of the torus (green).
Region labels as in Fig.~\ref{fig:bns}.}
\label{fig:ex3}
\end{figure}

The system's dynamics was studied, with TC-HAT \citep{Thota08b}, in the $(V,\lambda)$ plane, with the following values of the parameters: $m = 1022$, $K_s = 1e6$, $I = 678$, $A = 0.75$, $a=0.001$, $\mu_d = 1000$, $\alpha = 50$, $\mu_s = 1200$, $\Psi = 0.54219$, $\Phi = 0.60252$, $C = 6.5630e6$, $\mu_t = 1e5$, $r_0 = 0.4572$, $l=10$.  For large values of $V$, a grazing of a stable cycle with the boundary $x_2=0$ and a NS take place (blue and red in Fig.~\ref{fig:ex3}), and meet at
the border-NS (black) point.
Then, by systematically evaluating 1000 iterations (after transient) of the Poincar\'e map of the torus on a suitable cross-section, and by continuing the line on which the obtained torus image grazes the discontinuity boundary induced on the cross-section, we were able to trace an approximation of the grazing curve of the torus (green in Fig.~\ref{fig:ex3}).  More rigorous methods, based for example on discretisation of the invariant curve (see, e.g., \citep{Kevrekidis85,Dankowicz06}) could be used to obtain a more precise estimate of the quadratic coefficient.  This lies however beyond of the scope of this paper.  As predicted by the analysis carried out in Sect.~\ref{sec:bns}, the curve emanates from the codimension-two point tangentially to the NS curve. 

\section{Concluding remarks}
\label{sec:cr}
We have analysed the geometry of bifurcation curves around three codimension-two bifurcations in nonsmooth systems, namely the border-fold, the border-flip, and the border-Neimark-Sacker. Rather than aiming at the complete unfolding of the dynamics of a particular class of nonsmooth systems (e.g., piecewise smooth, impacting, or hybrid) dealing with a particular geometry of the involved discontinuity boundary (e.g., smooth or corner), we have focused on those results which are general to all scenarios.  Our approach applies to continuous-time as well as discrete-time systems, and basically consists of the analysis of a discrete-time (Poincar\'e) map defined only on one side of a boundary in its state space.  Explicit genericity conditions are listed and explained for each codimension-two case.

Of course, the weakness of this approach is that it cannot provide the complete unfolding of the bifurcation, but its power resides in its generality: as shown in the three examples that we have reported, it applies to a very broad class of nonsmooth systems and it may be relevant in various fields of science and engineering.

The natural sequel of this work would certainly aim at more detailed results, and possibly at the complete unfolding, of the codimension-two bifurcations analysed here, with specific reference to some smaller class of nonsmooth systems. 

\section*{Acknowledgements}
The first contributions on the topic of this paper were presented and discussed by Mario di Bernardo, Piotr Kowalczyk, and Yuri A. Kuznetsov in the context of piecewise smooth systems (informal meeting at the Bristol Centre for Applied Nonlinear Mathematics, University of Bristol, UK, summer 2003) and by Arne Nordmark for impacting system (at the meeting ``Piecewise smooth dynamical systems: Analysis, numerics and applications'', University of Bristol, Sept.~13--17, 2004).
The authors are grateful to M.~B., P.~K., Yu.~A.~K., and A.~N. for sharing their preliminary results, and to two anonymous reviewers whose criticisms significantly improved the paper.

\bibliographystyle{siam}
\bibliography{journalsshort,publishers,math,pws,discontinuous,ecology,mia_bib,mechanics}

\newpage

\appendix
\section{Border-fold bifurcation}
In the case of the border-fold bifurcation, conditions (i-iii) in Sec.~\ref{sec:fa}, expressed in the variable $u$ of the centre manifold, are summarised below:
\begin{list}{C}{}
 \item[(i.a)] $f_{uu}^0\neq 0$,
 \item[(i.b)] $f_{\alpha}^0\neq 0$,
 \item[(ii)] $h_u^0\neq 0$,
 \item[(iii)] $f_{uu}^0h_{\alpha_1}^0f_{\alpha_2}^0 - h_u^0f_{u\alpha_1}^0f_{\alpha_2}^0 \neq f_{uu}^0h_{\alpha_2}^0f_{\alpha_1}^0 - h_{u}^0f_{u\alpha_2}^0f_{\alpha_1}^0$
\end{list}
Note that (i.b) is redundant, since it is implied by (iii).

\subsection{Step one}
To reduce map (\ref{eq:bfd_cm}) to normal form we follow \citep{Kuznetsov04}, where however $\alpha\in\R$, while here $\alpha\in\R^2$.  The variable change $v=v(u,\alpha)$ is formally the same as in \citep{Kuznetsov04}, while  parameter change that we use is $\beta = \beta(\alpha)=|a(\mu(\alpha))|\mu(\alpha)$, $\mu_1(\alpha) = f_{0\alpha_1}^0\alpha_1+f_{0\alpha_2}^0\alpha_2+O(\|\alpha\|^2)$, $\mu_2(\alpha) = -f_{0\alpha_2}^0\alpha_1+f_{0\alpha_1}^0\alpha_2+O(\|\alpha\|^2)$, $a(\mu)=f_2(\alpha(\mu))+O(\|\alpha(\mu)\|)$,
with $a(0)=(1/2)f_{uu}^0\neq 0$ because of (i.a).  The inverse transformations have the following derivatives:
$$
u_{v}^0=\Frac{2}{|f_{uu}^0|},\quad
u_{\beta_2}^0=-\delta_{\alpha}^0\alpha_{\beta_2}^0,\quad
\delta_{\alpha}^0=\Frac{f_{u\alpha}^0}{f_{uu}^0},\quad
\alpha_{\beta_2}^0=\Frac{2}{|f_{uu}^0|\|f_{\alpha}^0\|^2}
\left[\begin{array}{c}
-f_{\alpha_2}^0\\
 f_{\alpha_1}^0
\end{array}\right].
$$

\subsection{Step two}
Consider the discontinuity boundary (\ref{eq:Sigma}).  The variable and parameter change $v=v(u,\alpha)$, $\beta=\beta(\alpha)$ is invertible near $(u,\alpha)=(0,0)$, so that condition (ii) implies that $h^{\mathrm{NF}}_v(0,0)=h_u^0 u_v^0\neq 0$, i.e., local existence and uniqueness, by the implicit function theorem, of a smooth function
$$
\sigma(\beta)=\sigma_{\beta_1}^0\beta_1+\sigma_{\beta_2}^0\beta_2+
O(\|\beta\|^2),
$$
such that $h^{\mathrm{NF}}(\sigma{\beta},\beta)=0$ for small $\|\beta\|$, so that the intersection of the discontinuity boundary $\mathcal{H}$ with the centre manifold $\mathcal{Z}^c$ is located at $v=\sigma(\beta)$.

We now prove, using condition (iii), that $\sigma_{\beta_2}^0\neq0$.  By differentiating both sides of $h^{\mathrm{NF}}(\sigma(\beta),\beta)=0$, i.e., of $h(u(\sigma(\beta),\beta),\alpha(\beta))=0$,
with respect to $\beta_2$, taking into account the derivatives in Appendix A.1, and evaluating at $\beta=0$ we get
$$
\sigma_{\beta_2}^0=
-\Frac{h_u^0 u_{\beta_2}^0+h_{\alpha}^0\alpha_{\beta_2}^0}
{h_u^0 u_{v}^0}=
\Frac{1}{h_u^0\|f_{\alpha}^0\|^2}\left(
\left(h_{\alpha_1}^0-\Frac{h_u^0f_{u\alpha_1}^0}{f_{uu}^0}\right)
f_{\alpha_2}^0-
\left(h_{\alpha_2}^0-\Frac{h_u^0f_{u\alpha_2}^0}{f_{uu}^0}\right)
f_{\alpha_1}^0
\right).
$$
Thanks to (i)--(iii), this ensures that $\sigma_{\beta_2}\neq 0$.

\subsection{Genericity conditions (ii) and (iii)}
In the original coordinates $z$ of map (\ref{eq:pm}), condition (ii) requires $H_z^0\nu^0\neq0$, where $\nu$ is the unit eigenvector of $F_z$ associated to the eigenvalue $1$.  

Consider now the fold curve defined by the system
\begin{equation}
\label{eq:gen_cond_fold_curve}
\begin{array}{rcl}
 F(z,\alpha)-z & = & 0,\\
 F_z(z,\alpha)\nu - \nu & = & 0,\\
 \langle \nu,\nu\rangle-1 & = & 0.
\end{array}
\end{equation}
In the space $(z,\nu,\alpha)$, condition (iii) means that the tangent vector to the fold curve is not tangent to the surface 
\begin{equation}
\label{eq:gen_cond_fold_plane}
H(z,\alpha)=0
\end{equation}
at $(z,\alpha)=(0,0)$.  The tangent vector to the fold curve is the null vector of the Jacobian of (\ref{eq:gen_cond_fold_curve}), so that bordering such Jacobian with the linearisation of (\ref{eq:gen_cond_fold_plane}) and imposing that the resulting square matrix is nonsingular at $(z,\nu,\alpha)=(0,\nu^0,0)$, i.e.,
$$
\det
\left(
\begin{array}{cccc}
F^0_z-I & 0 & F^0_{\alpha_1} & F^0_{\alpha_2}\\
F_{zz}^0\nu^0 & F^0_z-I &  F_{z\alpha_1}^0\nu^0 & F_{z\alpha_2}^0\nu^0\\
0 & 2(\nu^0)\ttt & 0 & 0\\
H^0_z & 0 & H^0_{\alpha_1} & H^0_{\alpha_2}
\end{array}
\right)\neq 0,
$$
we impose that the fold curve (\ref{eq:gen_cond_fold_curve}) intersects the surface (\ref{eq:gen_cond_fold_plane}) transversely, i.e., condition (iii).  This is nothing but requiring that the system (\ref{eq:gen_cond_fold_curve}), (\ref{eq:gen_cond_fold_plane}) be regular at $(z,\nu,\alpha)=(0,\nu^0,0)$. 

Equation (\ref{eq:gen_cond_fold_curve}), restricted to the centre manifold, becomes
$$
\begin{array}{rcl}
 f(u,\alpha)-u & = & 0,\\
 f_u(u,\alpha) - 1 & = & 0,
\end{array}
$$
and by the same reasoning, we obtain the condition 
$$
\det
\left(
\begin{array}{ccc}
f^0_u-1 & f^0_{\alpha_1} & f^0_{\alpha_2} \\
f^0_{uu} & f^0_{u\alpha_1} & f^0_{u\alpha_2}\\
h^0_u & h^0_{\alpha_1} & h^0_{\alpha_2}
\end{array}
\right)\neq 0,
$$
which is equivalent to (iii) since $f_u^0=1$ (fold condition).

\section{Border-flip bifurcation}
In the case of the border-flip bifurcation conditions (i-iii) in Sec.~\ref{sec:fa}, expressed in the variable $u$ of the centre manifold, are summarised below:
\begin{list}{C}{}
 \item[(i.a)] $\Frac{1}{2}(f_{uu}^0)^2+\Frac{1}{3}f_{uuu}^0\neq 0$,
 \item[(i.b)] $f_{u\alpha}^0\neq 0$,
 \item[(ii)] $h_u^0\neq 0$,
 \item[(iii)] $f_{u\alpha_1}^0h_{\alpha_2}^0\neq f_{u\alpha_2}^0h_{\alpha_1}^0$
\end{list}
Note that (i.b) is redundant, since it is implied by (iii).

\subsection{Step one}
Once again, to reduce map (\ref{eq:bfp_cm}) to normal form, we use the same variable change $v=v(u,\alpha)$ as in \citep{Kuznetsov04}, while the parameter change is $\beta_1 = \beta_1(\alpha) = g(\alpha_1,\alpha_2)$, $\beta_2 = \beta_2(\alpha) = h(0,\alpha)$, with $f_u(0,\alpha)=-(1+g(\alpha))$.  The inverse transformations have derivatives
$$
u_{v}^0=\Frac{1}{\sqrt{|c^0|}},\quad
u_{\beta_2}^0=0,\quad
\alpha_{\beta_2}^0=\Frac{1}{f^0_{u\alpha_1}h^0_{\alpha_2}-f^0_{u\alpha_2}h^0_{\alpha_1}}
\left[\begin{array}{c}
 -f_{u\alpha_2}^0\\
  f_{u\alpha_1}^0
\end{array}\right]
$$
with $c^0=(1/4)(f_{uu}^0)^2+(1/6)f_{uuu}^0\neq 0$ because of (i.a).

\subsection{Step two}
Consider the discontinuity boundary (\ref{eq:Sigma}).  As in the border-fold case, the variable and parameter change $v=v(u,\alpha)$, $\beta=\beta(\alpha)$ is invertible
near $(u,\alpha)=(0,0)$, so that condition (ii) implies
that $h^{\mathrm{NF}}_v(0,0)\neq 0$ and, by the implicit function theorem, that
the intersection of the discontinuity boundary $\mathcal{H}$
with the centre manifold $\mathcal{Z}^c$ is located at
$$
v=\sigma(\beta)=\sigma_{\beta_1}^0\beta_1+\sigma_{\beta_2}^0\beta_2+
O(\|\beta\|^2),
$$
for some smooth function $\sigma$.

The parameter change obviously makes $\sigma_{\beta_1}^0=0$.  We now prove that $\sigma_{\beta_2}^0\neq 0$.  By differentiating both sides of $h^{\mathrm{NF}}(\sigma(\beta),\beta)=0$,
i.e., of $h(u(\sigma(\beta),\beta),\alpha(\beta))=0$,
with respect to $\beta_2$, taking into account the derivatives in Appendix B.1 and evaluating at $\beta_2=0$ we get
$$
\sigma_{\beta_2}^0=
-\Frac{h_u^0u_{\beta_2}^0+h_{\alpha}^0\alpha_{\beta_2}^0}
{h_u^0 u_{v}^0}=
-\Frac{\sqrt{|c^0|}}{h_u^0},
$$
where condition (iii) ensures that $h_{\alpha}^0 \alpha_{\beta_2}=1$.  Thus (i)--(iii) imply that $\sigma_{\beta_2}^0\neq 0$.

\subsection{Genericity conditions (ii) and (iii)}
In the original coordinates $z$ of map (\ref{eq:pm}), condition (ii) requires $H_z^0\nu^0\neq0$, where $\nu$ is the unit eigenvector of $F_z$ associated to the eigenvalue $-1$. 

Consider now the flip curve defined by the system
\begin{equation}
\label{eq:gen_cond_flip_curve}
\begin{array}{rcl}
 F(z,\alpha)-z & = & 0,\\
 F_z(z,\alpha)\nu + \nu & = & 0,\\
 \langle \nu,\nu\rangle-1 & = & 0.
\end{array}
\end{equation}
Similarly to the border-fold case, condition (iii) is equivalent to
$$
\det
\left(
\begin{array}{cccc}
F^0_z-I & 0 & F^0_{\alpha_1} & F^0_{\alpha_2}\\
F_{zz}^0\nu^0 & F^0_z+I &  F_{z\alpha_1}^0\nu^0 & F_{z\alpha_2}^0\nu^0\\
0 & 2\nu\ttt & 0 & 0\\
H^0_z & 0 & H^0_{\alpha_1} & H^0_{\alpha_2}
\end{array}
\right)\neq 0.
$$
Equation (\ref{eq:gen_cond_flip_curve}), restricted to the centre manifold, becomes
$$
\begin{array}{rcl}
 f(u,\alpha)-u & = & 0,\\
 f_u(u,\alpha) + 1 & = & 0.
\end{array}
$$
Proceeding along the same lines we obtain the condition 
$$
\det
\left(
\begin{array}{ccc}
f^0_u-1 & f^0_{\alpha_1} & f^0_{\alpha_2} \\
f^0_{uu} & f^0_{u\alpha_1} & f^0_{u\alpha_2}\\
h^0_u & h^0_{\alpha_1} & h^0_{\alpha_2}
\end{array}
\right)\neq 0,
$$
which is equivalent to (iii) since $f^0_\alpha=0$ ($f(0,\alpha)=0$ by assumption) and $f^0_u=-1$ (flip condition).

\section{Border-Neimark-Sacker bifurcation}
In the case of the border-NS bifurcation, conditions (i-iii) in Sec.~\ref{sec:fa}, expressed in the variables $u$ of the centre manifold, are summarised below:
\begin{list}{C}{}
 \item[(i.a)] $e^{ik\theta^0}\neq 1$ for $k=1$, $2$, $3$, $4$,
 \item[(i.b)] the first Lyapunov coefficient of the NS normal form ($a^0$, see later) is nonzero,
 \item[(i.c)] $g_{\alpha}^0\neq 0$,
 \item[(ii)] $h_u^0\neq 0$,
 \item[(iii)] $g^0_{\alpha_1}h_{\alpha_2}^0 \neq g^0_{\alpha_2}h_{\alpha_1}^0$
\end{list}
Note that (i.c) is redundant, since it is implied by (iii).

\subsection{Step one}
Once again, to reduce map (\ref{eq:bns_cm}) to normal form, we use the same variable change $w=w(u,\alpha)$ (with $w=v_1+iv_2$) as in \citep{Kuznetsov04}, while the parameter change $\beta=\beta(\alpha)$ is formally the same as in B.1.
The inverse transformations $u=u(w,\bar{w},\beta)$ and $\alpha=\alpha(\beta)$ have derivatives
$$
u_{w}(0,0,0)=q^0,\quad u_{\bar{w}}(0,0,0)=\bar{q}^0,\quad
u_{\beta_2}(0,0,0)=0,\quad
\alpha_{\beta_2}^0=\Frac{1}{g^0_{\alpha_1}h^0_{\alpha_2}-g^0_{\alpha_2}h^0_{\alpha_1}}
\left[\begin{array}{c}
-g^0_{\alpha_2}\\
g^0_{\alpha_1}
\end{array}\right].
$$

\subsection{Step two}

\begin{figure}[t!]
\centerline{\includegraphics[width=1\textwidth]
{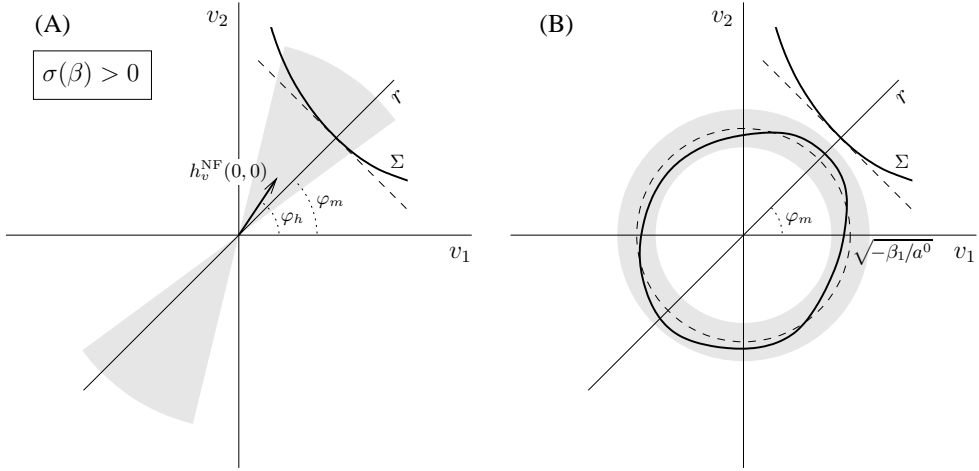}}
\caption{A. Local representation of the discontinuity boundary $\Sigma$ (thick line)
for small $\|v\|$ and $\|\beta\|$ as a straight (dashed) line tangent to
$\Sigma$ in the point of minimum distance of $\Sigma$ from the origin $v=0$
(case with $\sigma(\beta)>0$).
Since $\|\beta\|$ is small, the direction $\varphi_m$ of minimum distance is
close to the direction $\varphi_h$ of vector $h^{\mathrm{NF}}_v(0,0)$.
For $\varphi\in(\varphi_0,\varphi_1)$ (shaded area),
the discontinuity boundary $\Sigma$ can be represented in
coordinates $(r,\varphi)$.
B. The annular region (\ref{eq:bns_ar}) (shaded area) containing the
invariant curve (thick closed line) of the normal form map (\ref{eq:bns_nf})
and the (dashed) circle approached by the invariant curve
as $\beta \to 0$.}
\label{fig:bns_v}
\end{figure}

Denote by $\Sigma$ the discontinuity boundary (\ref{eq:Sigma}), where $v\in\R^2$.  Again, the variable and parameter change
$v=v(u,\alpha)$, $\beta=\beta(\alpha)$
that we used is invertible near $(u,\alpha)=(0,0)$, so that condition (ii)
implies that $h^{\mathrm{NF}}_v(0,0)=h_u^0 u_v^0\neq 0$,
where now $h^{\mathrm{NF}}_v(0,0)$ and $h_u^0$ are in $\R^2$ (row vectors)
and $u_v^0$ is a $2\times 2$ nonsingular matrix.
Geometrically, see Fig.~\ref{fig:bns_v}A, 
this means that
for small $\|v\|$ and $\|\beta\|$ we can represent the
discontinuity boundary (\ref{eq:Sigma}) as a straight line almost orthogonal
to $h^{\mathrm{NF}}_v(0,0)$ and slightly displaced from $v=0$ in the
direction of $h^{\mathrm{NF}}_v(0,0)$.

Let $\varphi_h$ be the angle of vector $h^{\mathrm{NF}}_v(0,0)$
with respect to axis $v_1$.
Technically,
$$
\varphi_h=\arctan_{2\pi}(h^{\mathrm{NF}}_{v1}(0,0),h^{\mathrm{NF}}_{v2}(0,0)),
$$
where $\arctan_{2\pi}$ is the four-quadrant inverse tangent in $[0,\,2\pi]$.
For any $\varphi$ in a neighbourhood of $\varphi_h$, introduce axis $r$
passing from the origin $v=0$ with direction $\varphi$,
so that positive and negative values of $r$ measure the distance
from the origin along directions $\varphi$ and $\varphi\pm\pi$,
respectively (see Fig.~\ref{fig:bns_v}A).
Coordinates $(r,\varphi)$ are like polar coordinates, but allow
differentiation with respect to $r$ at $r=0$.
We can therefore express the discontinuity boundary (\ref{eq:Sigma}) as
$$
\Sigma=\{(r,\varphi):
h^{\mathrm{NF}}((r\cos(\varphi),r\sin(\varphi)),\beta)=0\},
$$
where
$$
\left.
\Frac{d}{dr}h^{\mathrm{NF}}((r\cos(\varphi_h),r\sin(\varphi_h)),0)
\right|_{r=0}=
h^{\mathrm{NF}}_v(0,0)\left[
\begin{array}{c}\cos(\varphi_h)\\ \sin(\varphi_h)\end{array}\right]\neq 0
$$
(recall that, by definition of $\varphi_h$, $h^{\mathrm{NF}}_v(0,0)$ is
proportional to $(\cos(\varphi_h),\sin(\varphi_h))$), so that,
by the implicit function theorem, we can represent $\Sigma$ explicitly as
$r=\delta(\varphi,\beta)$, $\delta(\varphi,0)=0$,
for some smooth function $\delta$ defined for $\varphi$ in an
open neighbourhood $(\varphi_0,\varphi_1)$ of $\varphi_h$.

Now, define
$\varphi_m(\beta):=\mathrm{arg\,min}_{\varphi\in(\varphi_0,\varphi_1)}
\{|\delta(\varphi,\beta)|\}$ for $\beta\neq 0$ and note that
$\lim_{\beta\to 0}\varphi_m(\beta)=\varphi_h$, so that we can set
$\varphi_m^0=\varphi_h$.
Then, the minimum distance of $\Sigma$ from the origin $v=0$ is given by
the absolute value of
$$
\sigma(\beta):=\delta(\varphi_m(\beta),\beta)=
\sigma_{\beta_1}^0\beta_1+\sigma_{\beta_2}^0\beta_2+O(\|\beta\|^2),
$$
while its sign says whether the minimum is realised along the direction
$\varphi_m(\beta)$, if positive, or $\varphi_m(\beta)\pm\pi$, if negative.
In the first case (see Fig.~\ref{fig:bns_v}A), $v=0$ is a fixed point of the
normal form map (\ref{eq:bns_nf}), since $h^{\mathrm{NF}}(0,\beta)<0$,
while $v=0$ lies on the undescribed side of $\Sigma$ in the second case,
i.e., $h^{\mathrm{NF}}(0,\beta)>0$.

Similarly to the border-flip case, the parameter change implies that $\sigma_{\beta_1}^0=0$.  We now show that $\sigma_{\beta_2}^0\neq 0$.  By differentiating both sides of
$h^{\mathrm{NF}}((\delta(\varphi,\beta)\cos(\varphi),
\delta(\varphi,\beta)\sin(\varphi)),\beta)=0$, i.e., of
$$
h(u(\delta(\varphi,\beta)e^{i\varphi},
\delta(\varphi,\beta)e^{-i\varphi},\beta),\alpha(\beta))=0,
$$
with respect to $\beta_2$, taking into account the derivatives in Appendix C.1, and evaluating at $\beta_2=0$ we get
$$
\delta_{\beta_2}(\varphi,0)=
-\Frac{h_u^0u_{\beta_2}(0,0,0)+
h_{\alpha}^0\alpha_{\beta_2}^0}
{h_u^0(u_{w}^0e^{i\varphi}+u_{\bar{w}}^0e^{-i\varphi})}=
-\Frac{1}
{2h_u^0\mathrm{Re}(q^0e^{i\varphi})},
$$
which is well defined for $\varphi=\varphi_h$ thanks to (ii).
Indeed,
$u_{w}^0e^{i\varphi_h}+u_{\bar{w}}^0e^{-i\varphi_h}$ is nothing
but \linebreak $d/dr(u(re^{i\varphi_h},re^{-i\varphi_h},0))|_{r=0}$ and thus gives the
direction of $u$-perturbations from $u=0$ corresponding to $r$-perturbations
from $r=0$ along the direction $\varphi_h$, so that,
by definition of $\varphi_h$,
$\mathrm{Re}(q^0e^{i\varphi_h})$ is proportional to $h_u^0$.
Finally, we have
$$
\sigma_{\beta_2}^0=\delta_{\varphi}(\varphi_h,0)
\varphi_{m\beta_2}^0+\delta_{\beta_2}(\varphi_h,0)=
\delta_{\beta_2}(\varphi_h,0)
$$
(recall that $\delta(\varphi,0)=0$ for all $\varphi\in(\varphi_0,\varphi_1)$), so that $\sigma_{\beta_2}\neq0$ thanks to conditions (ii) and (iii) (which is necessary to show that $h_\alpha^0\alpha_{\beta_2}^0=1$).

Note that, in order to evaluate $\sigma_{\beta_2}^0$, we need an expression
for $\varphi_h$ in terms of variables $u$. For this we can write $u$ as a function of $(v,\beta)$, i.e.,
$$
u=u(v,\beta)=u(v_1+iv_2,v_1-iv_2,\beta)\\
$$
($u$ must be read as a function of $(w,\bar{w},\beta)$ in the right-most side),
so that
$$
\begin{array}{rcccl}
u_{v_1}^0 & = &
u_{w}(0,0,0)+u_{\bar{w}}(0,0,0) & = &
\phantom{-}2\mathrm{Re}(q^0),\\
u_{v_2}^0 & = &
u_{w}(0,0,0)i-u_{\bar{w}}(0,0,0)i & = & -2\mathrm{Im}(q^0),
\end{array}
$$
and
$$
\varphi_h=\arctan_{2\pi}\left(h_u^0u_{v_1}^0,h_u^0u_{v_2}^0\right)=
\arctan_{2\pi}\left(h_u^0\mathrm{Re}(q^0),-h_u^0\mathrm{Im}(q^0)\right).
$$

\subsection{Genericity conditions (ii) and (iii)}
Condition (ii) requires $H_z^0\left(\mathrm{Re}(nu^0),\mathrm{Im}(\nu^0)\right)\neq0$, where $\nu$ is the complex unit eigenvector of $F_z$ associated to the eigenvalue $(1+g)e^{i\theta}$.

The NS curve is described by the system 
\begin{equation}
\label{eq:gen_cond_ns_curve}
\begin{array}{rcl}
 F(z,\alpha)-z & = & 0,\\
 g(\alpha) & = & 0
\end{array}
\end{equation}
where, for any given $\alpha$, $g(\alpha)\in\R$ is obtained by solving the system
$$
\begin{array}{rcl}
 F_z(0,\alpha)\nu - (1+g)e^{i\theta}\nu & = & 0,\\
 \langle\nu,\nu\rangle -1 & = & 0,\\
 \mathrm{Re}(\nu)\ttt\mathrm{Im}(\nu) & = & 0,
\end{array}
$$
in the variables $(g,\theta,\nu)$.  In the space $(z,\alpha)$ condition (iii) means that the tangent vector to the NS curve is not tangent to the surface 
$$
H(z,\alpha)=0
$$
at $(z,\alpha)=(0,0)$.  Similarly to the border-fold and -flip cases, condition (iii) is equivalent to
$$
\det
\left(
\begin{array}{ccc}
F^0_z-I & F^0_{\alpha_1} & F^0_{\alpha_2}\\
g^0_z & g^0_{\alpha_1} & g^0_{\alpha_2}\\
H^0_z & H^0_{\alpha_1} & H^0_{\alpha_2}
\end{array}
\right)\neq 0.
$$
Equation (\ref{eq:gen_cond_ns_curve}), restricted to the centre manifold, becomes
$$
\begin{array}{rcl}
 f(u,\alpha)-u & = & 0,\\ g(\alpha) & = & 0.
\end{array}
$$
By the same reasoning we obtain the condition 
$$
\det
\left(
\begin{array}{ccc}
f^0_u-I & f^0_{\alpha_1} & f^0_{\alpha_2} \\
g^0_u & g^0_{\alpha_1} & g^0_{\alpha_2}\\
h^0_u & h^0_{\alpha_1} & h^0_{\alpha_2}
\end{array}
\right)\neq 0
$$
which is equivalent to (iii) since $f^0_\alpha=0$ ($f(0,\alpha)=0$ by assumption) and $f^0_u-I$ is nonsingular (condition (i.a) ($k=1$)). 

\subsection{Step three}
In this appendix we show that near $\beta=0$ the closed invariant curve of the normal form map (\ref{eq:bns_nf}) is contained in the parameter-dependent annular region (\ref{eq:bns_ar}) (we adapt the material from \citep{Kuznetsov10}, Chap. 5).

Assume the supercritical case, i.e., $a^0<0$, so that the
invariant curve exists for $\beta_1>0$ and is stable.
The annular region shrinks around the circle of equation
\begin{equation}
\label{eq:bns_cir}
\rho=\sqrt{-\Frac{\beta_1}{a(\beta)}},\quad \varphi\in[0,\,2\pi],
\end{equation}
with $O(\beta_1^\gamma)$-width (see Fig.~\ref{fig:bns_v}B) and map (\ref{eq:bns_nf_r}) maps $\rho$ into $\rho+\Delta\rho$ with
$\Delta\rho=\rho(\beta_1+a(\beta)\rho^2+\rho^3 R(\rho,\varphi,\beta))$ and
$$
\Delta\rho
\left\{\begin{array}{ll}
\ge\rho(2\beta_1^{\gamma+1/2}-\beta_1^{2\gamma}+O(\beta_1^{3/2})) &
\mathrm{if}\,0\le\rho\le
\sqrt{-\Frac{\beta_1}{a(\beta)}}(1-\beta_1^{\gamma-1/2}),\\[2mm]
\le\rho(-2\beta_1^{\gamma+1/2}-\beta_1^{2\gamma}+O(\beta_1^{3/2})) &
\mathrm{if}\,\rho\ge
\sqrt{-\Frac{\beta_1}{a(\beta)}}(1+\beta_1^{\gamma-1/2}).
\end{array}\right.
$$
Thus the orbits of map (\ref{eq:bns_nf}) enter the annular region if
$\gamma<1$ (the term $\beta_1^{\gamma+1/2}$ dominates the others and determines the sign of $\Delta\rho$), so that with $1/2<\gamma<1$ the stable invariant curve remains
in the annular region for small $\|\beta\|$.
Similarly, in the subcritical case, $a^0>0$, the
invariant curve exists for $\beta_1<0$ and is unstable, and the
orbits of map (\ref{eq:bns_nf}) exit the annular region if $\gamma<1$.
Again, with $1/2<\gamma<1$, the invariant curve remains
in the annular region for small $\|\beta\|$.
\end{document}